\documentclass[12pt]{amsart}
\usepackage{amssymb}

\theoremstyle{plain}
\newtheorem*{theorem}{Theorem}
\newtheorem{proposition}{Proposition}
\newtheorem{lemma}[proposition]{Lemma}

\theoremstyle{definition}
\newtheorem{example}{Example}
\newtheorem*{remark}{Remark}

\newcommand{\la}{\langle}
\newcommand{\ra}{\rangle}
\newcommand{\mat}{\begin{bmatrix}a&b\\c&d\end{bmatrix}}
\newcommand{\smat}{\left[\begin{smallmatrix}a&b\\c&d\end{smallmatrix}\right]}
\newcommand{\gie}{\mathbb{G}}
\renewcommand{\leq}{\leqslant}
\renewcommand{\geq}{\geqslant}

\begin{document}
\title[Frobenius theory]{Frobenius theory fails for semigroups \\
 of positive maps on von Neumann algebras}
\author[A. \L uczak]{Andrzej \L uczak}
\address{Faculty of Mathematics \\
        \L\'od\'z University \\
        ul. S. Banacha 22 \\
        90-238 \L\'od\'z, Poland}
\email{anluczak@math.uni.lodz.pl}
\thanks{Work supported by KBN grant 2PO3A 03024}
\date{}
\begin{abstract}
The eigenvectors of an ergodic semigroup of linear normal positive
unital maps on a von Neumann algebra are described. Moreover, it
is shown by means of examples, that mere positivity of the maps in
question is not sufficient for Frobenius theory as in \cite{A-H}
to hold.
\end{abstract}
\maketitle
\section{Introduction}\label{S0}
Frobenius theory for completely positive maps on von Neumann
algebras was developed in \cite{A-H} (further contributions to
this subject can be found in \cite{Gr} and \cite{Wa}). This theory
states, in paticular, that for an ergodic semigroup of completely
positive (or, in fact, even Schwarz) maps on a von Neumann algebra
its point spectrum forms a group, and the corresponding
eigenspaces are one-dimensional and spanned by a unitary operator.
The aim of this paper is to show that neither of these is true for a
semigroup of merely positive maps. Namely, we first prove that the
eigenvectors are either multiples of a partial isometry or linear
combinations of two partial isometries or multiples of a unitary
operator, and then show by means of examples, that there is a whole
class of semigroups of positive maps on a von Neumann algebra such
that their point spectra are not groups and the corresponding
eigenspaces have dimensions greater than one.
\section{Preliminaries and notation}\label{S1}
Let $M$ be a von Neumann algebra, and let $(\varPhi_g\colon
g\in\gie)$ be a semigroup of linear normal positive unital maps on
$M$. We shall be concerned with two cases: $\gie=\mathbb{N}_0$
--- all nonnegative integers, and $\gie=\mathbb{R}_+$ --- all
nonnegative reals (notice that in the first case the semigroup has
the form $(\varPhi^n\colon n=0,1,\dots)$, where $\varPhi$ is a
linear normal positive unital map on $M$).

A complex number $\lambda$ of modulus one is called an
\emph{eigenvalue} of the semigroup, if there is a nonzero $x\in M$
such that for each $g\in\gie$
\begin{equation}\label{e0}
 \varPhi_g(x) = \lambda^g x.
\end{equation}
The collection of all $x's$ such that \eqref{e0} holds is called
the \emph{eigenspace} corresponding to the eigenvalue $\lambda$,
and denoted by $M_{\lambda}$. In particular, $M_1$ is the
fixed-point space of the semigroup, and the semigroup is called
\emph{ergodic} if $M_1$ consists of multiples of the identity. The
set of all eigenvalues of the semigroup is called its \emph{point
spectrum}, and denoted by $\sigma((\varPhi_g))$. Let $\omega$ be a
normal faithful state on $M$ such that for each \linebreak
$g\in\gie,\ \omega\circ\varPhi_g = \omega$. The part of Frobenius
theory developed in \cite{A-H} which is of interest to us, states
that in this case, if we assume that \linebreak $(\varPhi_g\colon g\in\gie)$
is ergodic and the maps $\varPhi_g$ are two-positive, the
point spectrum is a group, and the corresponding eigenspaces are
one-dimensional and spanned by a unitary operator. A natural question
is if the same is true under the assumption of mere positivity of the
maps $\varPhi_g$. We shall show that this is not the case neither
for the group structure nor for the dimension of the eigenspaces.

Let $N$ be the $\sigma$-weak closure of the linear span of
$\bigcup_{\lambda}M_{\lambda}$, where the sum is taken over all
eigenvalues of $(\varPhi_g)$. Then according to \linebreak \cite[Theorem
1]{Lu} $N$ is a $JW^*$-algebra, by which is meant that $N$ is
\linebreak a $\sigma$-weakly closed linear space, closed with respect
to the Jordan product
\[
 x\circ y = \frac{1}{2}(xy + yx);
\]
moreover, $\varPhi_g|N$ are Jordan $^*$-automorphisms.
Consequently, if $x\in M_{\lambda}$, then
\[
 \varPhi_g(x^*) = \varPhi_g(x)^* = \bar{\lambda}^g x^*,
\]
meaning that $x^*\in M_{\bar{\lambda}}$, and for $x\in
M_{\lambda_1},\ y\in M_{\lambda_2}$ we have
\[
 \varPhi_g(x\circ y) = \varPhi_g(x)\circ\varPhi_g(y) = \lambda_1^g
 x\circ\lambda_2^g y = (\lambda_1\lambda_2)^g x\circ y,
\]
meaning that $x\circ y\in M_{\lambda_1\lambda_2}$. In particular,
for an eigenvector $x$ we have $x\circ x^*\in M_1$.
\section{Eigenvectors}\label{S2}
In the following theorem we shall describe the eigenvectors of
$(\varPhi_g)$.
\begin{theorem}
Assume that $(\varPhi_g)$ is ergodic, and let $x\in M_{\lambda}$
be an eigenvector of $(\varPhi_g)$. Then one
of the following possibilities holds:
\begin{enumerate}
\item[(i)] $x=\alpha v$, where $\alpha\in\mathbb{C}$, and $v$ is a
partial isometry in $M_{\lambda}$ such that
\[
 v^*v=e, \qquad vv^*=e^{\bot}
\]
for some nonzero projection $e$, $e\neq\boldsymbol{1}$;
\item[(ii)] $x=\alpha_1v_1+\alpha_2v_2$, where
$\alpha_1,\alpha_2\in\mathbb{C},\,\alpha_1\neq\alpha_2$, and
$v_1,v_2$ are partial isometries in $M_{\lambda}$ such that for
some nonzero projection $e,\ e\neq\boldsymbol{1}$,
\[
 v^*_1v_1=e, \quad v_1v^*_1=e^{\bot}, \quad v^*_2v_2=e^{\bot}, \quad
 v_2v^*_2=e;
\]
\item[(iii)] $x=\alpha u$, where $\alpha\in\mathbb{C}$, and $u$ is
a unitary operator in $M_{\lambda}$.
\end{enumerate}
\end{theorem}
\begin{proof}
We have $x\circ x^*\in M_1$, so renorming $x$, if
necessary, we may assume that
\begin{equation}\label{E1}
 x^*x+xx^*=\boldsymbol{1}.
\end{equation}
Multiplying both sides of the above equality by $x$ on the left
and on the right, respectively, we obtain
\begin{equation}\label{E2}
 xx^*x+x^2x^*=x=x^*x^2+xx^*x,
\end{equation}
so in particular $x^*$ commutes with $x^2$. Multiplying again both
sides of the second part of the above equality by $x^*$, we obtain
\begin{equation}\label{E3}
 x^*xx^*x+x^{*2}x^2=x^*x.
\end{equation}

(i) Assume that $x^2=0$. Then \eqref{E3} gives
\[
 (x^*x)^2=x^*x,
\]
so $x^*x$ is a projection $e$, and by \eqref{E1} $xx^*=e^{\bot}$,
thus $x$ is a partial isometry with initial projection $e$ and
final projection $e^{\bot}$, hence part (i) of the conclusion of
the theorem follows.

(ii) and (iii) Assume that $x^2\neq 0$. Then $x^2\in
M_{\lambda^{2}},\,x^{*2}\in M_{\bar{\lambda}^2}$, so $x^2\circ
x^{*2}\in M_1$, and since $x^2$ and $x^{*2}$ commute, we get
\[
 x^2x^{*2}=x^2\circ x^{*2}=\theta\boldsymbol{1}
\]
for some $\theta>0$, by the assumed ergodicity of $(\varPhi_g)$.
Denote $z=x^*x$. Then the above equality and \eqref{E3} give
\[
 z^2-z+\theta\boldsymbol{1}=0,
\]
that is
\[
 \left(z-\frac{1-\sqrt{1-4\theta}}{2}\,\boldsymbol{1}\right)
 \left(z-\frac{1+\sqrt{1-4\theta}}{2}\,\boldsymbol{1}\right)=0.
\]
It follows that $\theta\leq 1/4$, and we have two
possibilities for the spectrum of $z$: either $\theta=1/4$
in which case $\text{sp}\,z=\left\{1/2\right\}$, or $\theta
< 1/4$ in which case
$\text{sp}\,z=\left\{\frac{1-\sqrt{1-4\theta}}{2},
\,\frac{1+\sqrt{1-4\theta}}{2}\right\}$. The first possibility
gives at once
\[
 x^*x=\frac{1}{2}\,\boldsymbol{1}, \qquad
 xx^*=\frac{1}{2}\,\boldsymbol{1},
\]
and part (iii) of the conclusion of the theorem follows.

For the second possibility we have
\[
 x^*x=\frac{1-\sqrt{1-4\theta}}{2}\,e +
 \frac{1+\sqrt{1-4\theta}}{2}\,e^{\bot}
\]
where $e$ and $e^{\bot}$ are the spectral projections of $x^*x$.
From \eqref{E1} we get
\[
 xx^*=\frac{1+\sqrt{1-4\theta}}{2}\,e +
 \frac{1+\sqrt{1-4\theta}}{2}\,
 e^{\bot},
\]
and denoting
\[
 \alpha_1=\left(\frac{1-\sqrt{1-4\theta}}{2}\right)^{1/2}, \qquad
 \alpha_2=\left(\frac{1+\sqrt{1-4\theta}}{2}\right)^{1/2},
\]
we obtain
\begin{equation}\label{E4}
 \begin{aligned}
 |x|&=\alpha_1e + \alpha_2e^{\bot}\\
 |x^*|&=\alpha_2e + \alpha_1e^{\bot}.
 \end{aligned}
\end{equation}
Let
\[
 x=u|x|
\]
be the polar decomposition of $x$. Then
\[
 x^*=u^*|x^*|
\]
is the polar decomposition of $x^*$, and since $|x|$ and $|x^*|$ are
invertible the operator $u$ is unitary. Moreover,
\[
 |x^*|=u|x|u^*,
\]
which gives the equality
\[
 \alpha_2e + \alpha_1e^{\bot}=\alpha_1ueu^* + \alpha_2ue^{\bot}u^*,
\]
yielding
\[
 \alpha_1\boldsymbol{1} +
 (\alpha_2-\alpha_1)e=\alpha_2\boldsymbol{1} +
 (\alpha_1-\alpha_2)ueu^*.
\]
Thus
\[
 (\alpha_1-\alpha_2)\boldsymbol{1}=(\alpha_1-\alpha_2)(e + ueu^*),
\]
showing that
\[
 ueu^*=e^{\bot},
\]
and consequently,
\[
 e=\boldsymbol{1}-ueu^*=ue^{\bot}u^*.
\]
From the polar decomposition and formula \eqref{E4} we obtain
\begin{equation}\label{E4'}
 x=\alpha_1ue + \alpha_2ue^{\bot}=\alpha_1v_1 + \alpha_2v_2,
\end{equation}
where
\[
 v_1=ue, \qquad v_2=ue^{\bot}.
\]
We have
\[
 v^*_1v_1=e,\quad v_1v^*_1=ueu^*=e^{\bot},\quad
 v^*_2v_2=e^{\bot},\quad v_2v^*_2=ue^{\bot}u^*=e,
\]
so \eqref{E4'} is the representation of $x$ as in part (ii) of the
conclusion of the theorem. It remains to show that $v_1,v_2\in
M_{\lambda}$. Equality \eqref{E2} gives
\[
 xx^*x=x-x^2x^*=x-x^2\circ x^*
\]
because $x^2$ and $x^*$ commute. Since $x^2\in M_{{\lambda}^2}$ and
$x^*\in M_{\bar{\lambda}}$, we get $x^2\circ x^*\in M_{\lambda}$, so
$xx^*x\in M_{\lambda}$. We have
\[
 xx^*x=(\alpha_1v_1 + \alpha_2v_2)(\alpha^2_1e + \alpha^2_2e^{\bot})
 =\alpha^3_1v_1 + \alpha^3_2v_2,
\]
and the equality
\[
 \varPhi_g(xx^*x)=\lambda^g xx^*x
\]
yields
\begin{equation}\label{E5}
 \alpha^3_1\varPhi_g(v_1) +\alpha^3_2\varPhi_g(v_2)
 =\lambda^g\alpha^3_1v_1 + \lambda^g\alpha^3_2v_2.
\end{equation}
On the other hand $\varPhi_g(x)=\lambda^g x$, which gives
\begin{equation}\label{E6}
 \alpha_1\varPhi_g(v_1) + \alpha_2\varPhi_g(v_2) = \lambda^g\alpha_1v_1
 + \lambda^g\alpha_2v_2.
\end{equation}
Multiplying both sides of equality \eqref{E6} by $\alpha^2_2$ and
substracting \eqref{E6} from \eqref{E5} we obtain
\[
 \alpha_1(\alpha^2_1-\alpha^2_2)\varPhi_g(v_1) =
 \lambda^g\alpha_1(\alpha^2_1-\alpha^2_2)v_1,
\]
which gives $\varPhi_g(v_1)=\lambda^g v_1$,
and analogously $\varPhi_g(v_2)=\lambda^g v_2$.
\end{proof}
\section{Positivity of elements from $\textbf{Mat}_2(M)$}\label{S3}
In what follows we shall need a number of properties concerning
positivity of matrices with elements in a von Neumann algebra.
They will be exploited in examples in a particular case of an
abelian von Neumann algebra, but as these properties seem to be
interesting in their own right, we prove them here in slightly
greater generality. For more information on positivity of such
matrices the reader is referred to \cite[Chapter IV.3]{Ta} and
\cite{Wo}. It should be added that some of the facts obtained
below can be given alternative proofs based on methods used in
\cite{Ta,Wo}.

Let $M$ be a von Neumann algebra, and let
$\widetilde{M}=\textbf{Mat}_2(M)$ be the algebra of
$2\times2$-matrices with elements from $M$. Assuming that $M$ acts
on a Hilbert space $\mathcal{H}$, we can consider $\widetilde{M}$
as acting on the Hilbert space $\widetilde{\mathcal{H}} =
\mathcal{H}\oplus\mathcal{H}$.

\begin{proposition}\label{P1}
Let $A=\left[\begin{smallmatrix}a & b\\ c &
d\end{smallmatrix}\right]\in\widetilde{M}$.
\begin{enumerate}
\item[(i)]$A\geq 0$ if and only if $a,d\geq 0,\,c=b^*$, and for
each $\varepsilon > 0 \linebreak d\geq
b^*(a+\varepsilon\boldsymbol{1})^{-1}b$.
\item[(ii)] Assume that $a,b\text{ and }c$ commute. Then
$A\geq 0$ if and only if $a,d\geq 0,\,c=b^*$, and $ad\geq b^*b$.
\end{enumerate}
\end{proposition}
\begin{proof}
Calculate first the quadratic form of $A$. For
$\tilde{\xi}=\left(\begin{smallmatrix}\xi_1 \\
\xi_2\end{smallmatrix}\right)$,\linebreak
$\xi_1,\xi_2\in\mathcal{H}$, we have
\begin{equation}\label{e1}
 \begin{aligned}
 \la
 A\tilde{\xi},\tilde{\xi}\ra_{\widetilde{\mathcal{H}}}&=\la\begin{pmatrix}
 a\xi_1+b\xi_2 \\ c\xi_1+d\xi_2\end{pmatrix},\begin{pmatrix}\xi_1
 \\ \xi_2\end{pmatrix}\ra_{\widetilde{\mathcal{H}}} \\
 &=\la a\xi_1,\xi_1\ra_{\mathcal{H}}+\la
 b\xi_2,\xi_1\ra_{\mathcal{H}}+\la
 c\xi_1,\xi_2\ra_{\mathcal{H}}+\la d\xi_2,\xi_2\ra_{\mathcal{H}}.
 \end{aligned}
\end{equation}

It is clear that in order that $A$ be positive we must have
$a,d\geq 0$ and $c=b^*$, so we restrict attention to matrices $A$
of the form $A=\left[\begin{smallmatrix}a & b \\ b^* & d
\end{smallmatrix}\right]$ with $a,d \geq 0$. Then \eqref{e1}
becomes
\begin{equation}\label{e2}
 \la A\tilde{\xi},\tilde{\xi}\ra_{\widetilde{\mathcal{H}}}=\la
 a\xi_1,\xi_1\ra_{\mathcal{H}}+\la b\xi_2,\xi_1\ra_{\mathcal{H}}
 +\la\xi_1,b\xi_2\ra_{\mathcal{H}} + \la
 d\xi_2,\xi_2\ra_{\mathcal{H}}.
\end{equation}

(i) \emph{Step 1}. First we shall show that for $A$ of the form
$A=\left[\begin{smallmatrix}\boldsymbol{1} & b \\ b^* & z
\end{smallmatrix}\right]$, $A\geq 0$ if and only if $z\geq b^*b$.
This is virtually proved in\\ \cite[Lemma 3.1]{Wo}. For the sake of
completeness we give a simple proof below.

For the quadratic form of $A$ we have
\begin{equation}\label{e3}
 \la
 A\tilde{\xi},\tilde{\xi}\ra_{\widetilde{\mathcal{H}}}=
 \la\xi_1,\xi_1\ra_{\mathcal{H}}+\la
 b\xi_2,\xi_1\ra_{\mathcal{H}}+\la\xi_1,b\xi_2\ra_{\mathcal{H}}+
 \la z\xi_2,\xi_2\ra_{\mathcal{H}}.
\end{equation}

Let $b^*b\leq z$. Then
\begin{align*}
 \la A\tilde{\xi},\tilde{\xi}\ra_{\widetilde{\mathcal{H}}}&\geq
 \la\xi_1,\xi_1\ra_{\mathcal{H}}+\la b\xi_2,\xi_1\ra_{\mathcal{H}}
 +\la\xi_1,b\xi_2\ra_{\mathcal{H}}+\la
 b\xi_2,b\xi_2\ra_{\mathcal{H}}\\&=\la\xi_1+b\xi_2,\xi_1+b\xi_2\ra_{\mathcal{H}}
 \geq 0,
\end{align*}
showing that $A\geq 0$.

Conversely, if $A\geq 0$ then substituting $-\xi_1$ for $\xi_1$ in
\eqref{e3}, we obtain
\[
 0\leq \la\xi_1,\xi_1\ra_{\mathcal{H}}-\la
 b\xi_2,\xi_1\ra_{\mathcal{H}}
 -\la\xi_1,b\xi_2\ra_{\mathcal{H}}+\la
 z\xi_2,\xi_2\ra_{\mathcal{H}},
\]
that is
\[
\la b\xi_2,\xi_1\ra_{\mathcal{H}}+\la\xi_1,b\xi_2\ra_{\mathcal{H}}
-\la \xi_1,\xi_1\ra_{\mathcal{H}}\leq\la
z\xi_2,\xi_2\ra_{\mathcal{H}}.
\]

Now putting $\xi_1=b\xi_2$ in the above inequality, we get
\[
 \la b\xi_2,b\xi_2\ra_{\mathcal{H}}\leq\la
 z\xi_2,\xi_2\ra_{\mathcal{H}},
\]
which shows that $b^*b\leq z$.

\emph{Step 2}. Let now $A=\left[\begin{smallmatrix}a & b \\ b^* &
d\end{smallmatrix}\right]$ be arbitrary. $A\geq 0$ if and only if
for each $\varepsilon > 0, \quad
A_{\varepsilon}=\left[\begin{smallmatrix}a+\varepsilon\boldsymbol{1}&
b \\ b^*& d \end{smallmatrix}\right]\geq 0$, and denoting
$a_\varepsilon = a + \varepsilon\boldsymbol{1}$, we obtain
\[
 \la
 A_{\varepsilon}\tilde{\xi},\tilde{\xi}\ra_{\widetilde{\mathcal{H}}}=
 \la a_{\varepsilon}\xi_1,\xi_1\ra_{\mathcal{H}}+\la b\xi_2,\xi_1
 \ra_{\mathcal{H}}+\la\xi_1,b\xi_2\ra_{\mathcal{H}}+\la
 d\xi_2,\xi_2\ra_{\mathcal{H}}.
\]

Putting $\eta_1=a_{\varepsilon}^{1/2}\xi_1$, we get
\[
 \begin{aligned}
 \la
 A_{\varepsilon}\tilde{\xi},\tilde{\xi}\ra_{\widetilde{\mathcal{H}}}&=
 \la \eta_1,\eta_1\ra_{\mathcal{H}}+\la
 a_{\varepsilon}^{-1/2}b\xi_2,\eta_1\ra_{\mathcal{H}}
 +\la\eta_1,a_{\varepsilon}^{-1/2}b\xi_2\ra_{\mathcal{H}}+\la
 d\xi_2,\xi_2\ra_{\mathcal{H}}\\
 &=\la\begin{bmatrix}\boldsymbol{1}& a_{\varepsilon}^{-1/2}b \\
 b^* a_{\varepsilon}^{-1/2} & d \end{bmatrix}\begin{pmatrix}\eta_1\\
 \xi_2\end{pmatrix},\begin{pmatrix}\eta_1 \\
 \xi_2\end{pmatrix}\ra_{\widetilde{\mathcal{H}}}.
 \end{aligned}
\]

Since $a_{\varepsilon}^{1/2}$ maps $\mathcal{H}$ onto
$\mathcal{H}$ in a 1--1 way, we see that $A_{\varepsilon}\geq 0$
if and only if $\left[\begin{smallmatrix}\boldsymbol{1} &
a_{\varepsilon}^{-1/2}b \\ b^* a_{\varepsilon}^{-1/2} & d
\end{smallmatrix}\right]\geq 0$, which by \emph{Step 1} is
equivalent to the condition
\[
 d \geq
 b^*a_{\varepsilon}^{-1}b=b^*(a+\varepsilon\boldsymbol{1})^{-1}b,
\]
and the proof of (i) is complete.

(ii) The reasoning is similar to that in part (i) using the simple
\linebreak $\varepsilon$-trick. Namely, the inequality $ad\geq
b^*b$ is equivalent to $(a+\varepsilon\boldsymbol{1})d \geq b^*b$
for each $\varepsilon >0$, which in turn, by the assumed
commutation property, is equivalent to $d\geq
b^*(a+\varepsilon\boldsymbol{1})^{-1}b$. Applying part (i)
finishes the proof.
\end{proof}
\begin{remark}
In virtually the same way we obtain the following variant of (i):
\begin{enumerate}
\item[(i${}'$)] $A\geq 0$ if and only if $a,d\geq 0,\quad c=b^*$,
and for each $\varepsilon > 0 \linebreak a\geq
b(d+\varepsilon\boldsymbol{1})^{-1}b^*$.
\end{enumerate}
\end{remark}
\begin{lemma}\label{L2}
Let $\left[\begin{smallmatrix}a & b \\ b^* & d
\end{smallmatrix}\right]\geq 0$. Then
\begin{enumerate}
\item[(i)]$\left[\begin{smallmatrix}d & b^* \\ b & a
\end{smallmatrix}\right]\geq 0$.
\item[(ii)] For each $x,y\in M \quad
\left[\begin{smallmatrix}xax^* & xby^* \\ yb^*x^* & ydy^*
\end{smallmatrix}\right]\geq 0$.
\end{enumerate}
\end{lemma}
\begin{proof}
(i) follows from the equality
\[
 \begin{bmatrix}d & b^* \\ b & a \end{bmatrix} = \begin{bmatrix}0
 & \boldsymbol{1} \\ \boldsymbol{1} & 0
 \end{bmatrix}\begin{bmatrix}a & b \\ b^* & d \end{bmatrix}
 \begin{bmatrix}0 & \boldsymbol{1} \\ \boldsymbol{1} & 0
 \end{bmatrix},
\]
and (ii) from the equality
\[
 \begin{bmatrix}xax^* & xby^* \\ yb^*x^* & ydy^* \end{bmatrix}=
 \begin{bmatrix}x & 0 \\ 0 & y \end{bmatrix}\begin{bmatrix}a & b
 \\ b^* & d \end{bmatrix}\begin{bmatrix}x^* & 0 \\ 0 & y^*
 \end{bmatrix}.
\]
\end{proof}
\begin{lemma}\label{L3}
Let $a$ commute with $b$, and assume further that either $b$ is
normal or that $b$ commutes with $d$. If
$\left[\begin{smallmatrix}a & b \\ b^* & d
\end{smallmatrix}\right]\geq 0$, then $\left[\begin{smallmatrix}a
& b^* \\ b & d \end{smallmatrix}\right]\geq 0$.
\end{lemma}
\begin{proof}
Let $\left[\begin{smallmatrix}a & b \\ b^* & d
\end{smallmatrix}\right]\geq 0$, and assume first that $b$ is
normal. Then by Proposition \ref{P1} (i) we have for each
$\varepsilon > 0$
\[
 d \geq b^*(a+\varepsilon\boldsymbol{1})^{-1}b =
 b^*b(a+\varepsilon\boldsymbol{1})^{-1} =
 b(a+\varepsilon\boldsymbol{1})^{-1}b^*,
\]
which again by Proposition \ref{P1} (i) means that
$\left[\begin{smallmatrix}a & b^* \\ b & d
\end{smallmatrix}\right]\geq 0$.

Now let $b$ commute with $d$. Then for each $\varepsilon > 0$,
\[
 d \geq b^*(a+\varepsilon\boldsymbol{1})^{-1}b =
 |b|^2(a+\varepsilon\boldsymbol{1})^{-1} =
 |b|(a+\varepsilon\boldsymbol{1})^{-1}|b|,
\]
so $\left[\begin{smallmatrix}a & |b| \\ |b| & d
\end{smallmatrix}\right]\geq 0$. By Lemma \ref{L2} (i) it follows
that $\left[\begin{smallmatrix}d & |b| \\ |b| & a
\end{smallmatrix}\right]\geq 0$, and thus, on account of
Proposition \ref{P1} (i), for each $\varepsilon > 0$,
\[
 a \geq |b|(d+\varepsilon\boldsymbol{1})^{-1}|b| =
 b^*b(d+\varepsilon\boldsymbol{1})^{-1} =
 b^*(d+\varepsilon\boldsymbol{1})^{-1}b,
\]
which means that $\left[\begin{smallmatrix}d & b \\ b^* & a
\end{smallmatrix}\right]\geq 0$. Applying again Lemma \ref{L2} (i)
we obtain $\left[\begin{smallmatrix}a & b^* \\ b & d
\end{smallmatrix}\right]\geq 0$.
\end{proof}
\section{Examples}\label{S4}
Let us begin with a simple example.
\begin{example}\label{Ex1}
Keeping the notation of Section \ref{S3}, put
$\widetilde{M}=\mathbb{B}(\mathbb{C}^2)$ (i.e. $M=\mathbb{C}$),
$\omega=\frac{1}{2}tr$, and let $\lambda_0\in\mathbb{C}$ be such
that $|\lambda_0|=1,\ \lambda_0\ne1$. Define
$\varPhi:\widetilde{M}\to\widetilde{M}$ as
\[
 \varPhi\left(\mat\right) =
 \begin{bmatrix}\frac{a+d}{2}&\lambda_0b\\
 \bar{\lambda}_0c&\frac{a+d}{2}\end{bmatrix}.
\]
It is clear that $\varPhi$ is linear normal unital, and that
$\omega\circ\varPhi=\omega$. Moreover, for $\smat\geq 0$, we have
$a,d \geq 0,\ c=\bar{b}$, and $\left[\begin{smallmatrix}d & b \\
\bar{b} & a \end{smallmatrix}\right]\geq 0$, thus \linebreak
$\left[\begin{smallmatrix}a+d & 2b \\ 2\bar{b} & a+d
\end{smallmatrix}\right]\geq 0$, so
$\left[\begin{smallmatrix}\frac{a+d}{2} & b \\ \bar{b} &
\frac{a+d}{2} \end{smallmatrix}\right]\geq 0$, and consequently
$\left[\begin{smallmatrix}\frac{a+d}{2} & \lambda_0b \\
\bar{\lambda}_0\bar{b} & \frac{a+d}{2}
\end{smallmatrix}\right]\geq 0$, showing that $\varPhi$ is
positive.

The equality
\[
 \varPhi\left(\mat\right) = \mat
\]
yields
\[
 \frac{a+d}{2} = a = d,\qquad \lambda_0b = b,\qquad
 \bar{\lambda}_0c = c,
\]
hence $b=c=0$, and the fixed-points have the form
$\left[\begin{smallmatrix}a & 0 \\ 0 & a \end{smallmatrix}\right]
= a\left[\begin{smallmatrix}1 & 0 \\ 0 & 1
\end{smallmatrix}\right]$, which means that $(\varPhi^n)$ is
ergodic.

(i) Let $\lambda_0$ be such that $\lambda_0\ne-1,\
\lambda_0^3\ne1$, and let $\lambda\ne1$ be an eigenvalue of
$(\varPhi^n)$. Then
\[
 \frac{a+d}{2} = \lambda a = \lambda d,\qquad \lambda_0b = \lambda
 b,\qquad \bar{\lambda}_0c = \lambda c,
\]
which yields
\[
 a = d = 0,\qquad \bar{\lambda}\lambda_0b = b,\qquad
 \bar{\lambda}\bar{\lambda}_0c = c.
\]
Thus either $\lambda = \lambda_0,\ c = 0$ or $\lambda =
\bar{\lambda}_0,\ b = 0$, so $\lambda_0 \text{ and }
\bar{\lambda}_0$ are the only eigenvalues of $(\varPhi^n)$
different from $1$, with the eigenspaces
\[
 \widetilde{M}_{\lambda_0} = \left\{\begin{bmatrix}0 & b \\ 0 & 0
 \end{bmatrix}\colon b\in\mathbb{C}\right\}, \qquad
 \widetilde{M}_{\bar{\lambda}_0} = \left\{\begin{bmatrix}0 & 0 \\
 c & 0 \end{bmatrix}\colon c\in\mathbb{C}\right\}.
\]
Consequently, $\sigma((\varPhi^n)) =
\{1,\lambda_0,\bar{\lambda}_0\}$, which is not a group if
\linebreak $\lambda_0\ne-1,\ \lambda_0^3\ne1$.

(ii) Now let $\lambda_0 = -1$. The above calculations give
$\sigma((\varPhi^n)) = \{1,-1\}$, and
\[
 \widetilde{M}_{-1} = \left\{\begin{bmatrix} 0 & b \\ c & 0
 \end{bmatrix}\colon b,c\in\mathbb{C}\right\},
\]
so the eigenspace is not one-dimensional.\hfill\qedsymbol
\end{example}

Let us observe that the reasoning above may be repeated with
virtually no change for the semigroup $(\varPhi_t\colon t\geq 0)$
defined as
\[
 \varPhi_t\left(\mat\right) = \begin{bmatrix} \frac{a+d}{2} &
 \lambda_0^t b \\ \bar{\lambda}_0^t c & \frac{a+d}{2}
 \end{bmatrix},
\]
thus giving a corresponding example in the continuous case.
\begin{remark}
The triple $(\mathbb{B}(\mathbb{C}^2),(\varPhi_t),\omega)$ from
the above example constitutes what in \cite{Gr} is called an
irreducible $W^*$-dynamical system. In \linebreak\cite[Theorem
3.8]{Gr} it is proved that under the assumption that the
$\varPhi_t$'s are Schwarz maps every such a system on a full
algebra has trivial point spectrum (i.e. consisting only of $1$).
As we see this is not the case if we assume only positivity of the
maps $\varPhi_t$'s.
\end{remark}
Now we construct a more involved
example (in fact, a class of examples) in which we shall see that
all the possibilities for the eigenvectors given in Theorem may
occur.
\begin{example}\label{Ex2}
Let $M$ be abelian, let $\omega$ be a normal faithful state on
$M$, and let $\varPsi$ be a positive normal unital map on $M$ such
that $\omega\circ\varPsi = \omega,\ (\varPsi^n)$ is ergodic, and
$\sigma((\varPsi^n)) = \{-1,1\}$. The abelianess of $M$ implies
that $\varPsi$ is completely positive (cf. \cite[Chapter IV.3]{Ta}),
thus according to \cite{A-H} the eigenspace corresponding to $-1$
is one-dimensional and spanned by a unitary operator $u$, i.e.
\[
 \varPsi(x) = -x
\]
if and only if $x$ is a multiple of $u$.

Put $\widetilde{M} = \textbf{Mat}_2(M)$,
\[
 \tilde\omega\left(\mat\right) = \frac{1}{2}[\omega(a) +
 \omega(d)],
\]
and let $\lambda_0\in\mathbb{C}$ be such that $|\lambda_0|=1,\
\lambda_0\notin\{-1,1\}$. Define
$\varPhi:\widetilde{M}\to\widetilde{M}$ as
\[
 \varPhi\left(\mat\right) =
 \begin{bmatrix}\varPsi\left(\frac{a+d}{2}\right) & \lambda_0\varPsi(b) \\
 \bar{\lambda}_0\varPsi(c) & \varPsi\left(\frac{a+d}{2}\right)\end{bmatrix}.
\]
$\varPhi$ is a linear normal unital map on $\widetilde{M}$ such
that $\tilde\omega\circ\varPhi = \tilde\omega$. Arguing as in
Example \ref{Ex1}, and using Lemma \ref{L2} with $x =
\lambda_0\boldsymbol{1},\ y = \boldsymbol{1}$, and Lemma \ref{L3}, we see
that if $\smat\geq 0$, then
$\left[\begin{smallmatrix}\frac{a+d}{2} & \lambda_0 b \\ \bar{\lambda}_0 c &
\frac{a+d}{2}\end{smallmatrix}\right]\geq 0$, so $\varPhi$ is
positive by virtue of the complete positivity of $\varPsi$.

The equality
\[
 \varPhi\left(\mat\right) = \mat
\]
yields
\[
 \varPsi\left(\frac{a+d}{2}\right) = a = d,\qquad \lambda_0\varPsi(b) =
 b,\qquad \bar{\lambda}_0\varPsi(c) = c,
\]
hence $b = c = 0$ and $\varPsi(a) = a$. From the ergodicity of
$\varPsi$ it follows that $a$ is a multiple of $\boldsymbol{1}$, so
the fixed-points of $(\varPhi^n)$ have the form
$\left[\begin{smallmatrix}\theta\boldsymbol{1} & 0 \\ 0 &
\theta\boldsymbol{1}\end{smallmatrix}\right]$ with
$\theta\in\mathbb{C}$, which means that $(\varPhi^n)$ is ergodic.
Let $\lambda\ne1$ be an eigenvalue of $(\varPhi^n)$. Then
\[
 \varPsi\left(\frac{a+d}{2}\right) = \lambda a = \lambda d,\qquad
 \lambda_0\varPsi(b) = \lambda b,\qquad \bar{\lambda}_0\varPsi(c)
 = \lambda c,
\]
which yields
\begin{equation}\label{e4}
 \varPsi(a) = \lambda a,\qquad \varPsi(b) = \bar{\lambda}_0\lambda b,
 \qquad \varPsi(c) = \bar{\lambda}_0\lambda c.
\end{equation}

(i) Take $\lambda_0$ such that $\lambda_0\neq i,\
\lambda_0\neq-i,\ \lambda^3_0\neq1,\ \lambda^3_0\neq-1$.
Equalities \eqref{e4} yield the following possibilities:

(i.1) $\lambda=-1$. Then $a=d$ is a multiple of $u,\
b=c=0$, and
\[
 \widetilde{M}_{-1}=\left\{\alpha\begin{bmatrix}u & 0\\ 0 & u
 \end{bmatrix}\colon \alpha\in\mathbb{C}\right\},
\]
so the eigenvector corresponding to the eigenvalue $-1$ is as in
part (iii) of Theorem.

(i.2) $\lambda\neq-1$. Then $a=d=0$, and one of the four
situations must occur:

(i.2.1) $\lambda=\lambda_0$. Then $b$ is a multiple of
$\boldsymbol{1},\ c=0$, and
\[
 \widetilde{M}_{\lambda_0}=\left\{\alpha\begin{bmatrix}0 &
 \boldsymbol{1}\\ 0 & 0
 \end{bmatrix}\colon\alpha\in\mathbb{C}\right\},
\]
so the eigenvector corresponding to the eigenvalue $\lambda_0$ is
as in part (i) of Theorem.

(i.2.2) $\lambda=-\lambda_0$. Then $b$ is a multiple of $u,\ c=0$, and
\[
 \widetilde{M}_{-\lambda_0}=\left\{\alpha\begin{bmatrix}0 & u\\
 0 & 0\end{bmatrix}\colon\alpha\in\mathbb{C}\right\}.
\]

(i.2.3) $\lambda=\bar{\lambda}_0$. Then $b=0,\ c$ is a multiple of
$\boldsymbol{1}$, and
\[
 \widetilde{M}_{\bar{\lambda}_0}=\left\{\alpha\begin{bmatrix}0 &
 0\\ \boldsymbol{1} & 0
 \end{bmatrix}\colon\alpha\in\mathbb{C}\right\}.
\]

(i.2.4) $\lambda=-\bar{\lambda}_0$. Then $b=0,\ c$ is a multiple
of $u$, and
\[
 \widetilde{M}_{-\bar{\lambda}_0}=\left\{\alpha\begin{bmatrix}0 &
 0\\u & 0 \end{bmatrix}\colon\alpha\in\mathbb{C}\right\}.
\]

Moreover, we have
\[
 \sigma((\varPhi^n))=\{1,-1,\lambda_0,\bar{\lambda}_0,-\lambda_0,
 -\bar{\lambda}_0\},
\]
which is not a group under our assumptions on $\lambda_0$.

Now take $\lambda_0=i$. Equations \eqref{e4} become then
\[
 \varPsi(a)=\lambda a,\qquad \varPsi(b)=-i\lambda b,\qquad
 \varPsi(c)=i\lambda c.
\]
As in part (i) we have the possibilities:

(ii.1) $\lambda=-1$, in which case $a=d$ is a multiple of
$\boldsymbol{1}$, and $b=c=0$.

(ii.2) $\lambda\neq-1$, in which case $a=d=0$, and we may only
have either $\lambda=i$ or $\lambda=-i$. In the first case $b$ is
a multiple of $\boldsymbol{1},\ c$ is a multiple of $u$, and
\[
 \widetilde{M}_i=\left\{\alpha_1\begin{bmatrix}0 &
 \boldsymbol{1}\\ 0 & 0\end{bmatrix}+ \alpha_2\begin{bmatrix}0 &
 0\\ u &
 0\end{bmatrix}\colon\alpha_1,\alpha_2\in\mathbb{C}\right\},
\]
so the situation is as in part (ii) of Theorem with
\[
 v_1=\begin{bmatrix}0 & \boldsymbol{1}\\ 0 & 0\end{bmatrix},\
 v_2=\begin{bmatrix}0 & 0\\ u & 0\end{bmatrix},\
 e=\begin{bmatrix}0 & 0\\ 0 & \boldsymbol{1}\end{bmatrix},\
 e^{\bot}=\begin{bmatrix}\boldsymbol{1} & 0\\ 0 & 0\end{bmatrix}.
\]
In the second case $b$ is a multiple of $u,\ c$ is a multiple of
$\boldsymbol{1}$, and
\[
 \widetilde{M}_{-i}=\left\{\alpha_1\begin{bmatrix}0 & u\\ 0 &
 0\end{bmatrix}+ \alpha_2\begin{bmatrix}0 & 0\\ \boldsymbol{1} &
 0\end{bmatrix}\colon\alpha_1,\alpha_2\in\mathbb{C}\right\},
\]
so again part (ii) of Theorem occurs with
\[
 v_1=\begin{bmatrix}0 & u\\ 0 & 0\end{bmatrix},\
 v_2=\begin{bmatrix}0 & 0\\ \boldsymbol{1} & 0\end{bmatrix},\
 e=\begin{bmatrix}\boldsymbol{1} & 0\\ 0 & 0\end{bmatrix},\
 e^{\bot}=\begin{bmatrix}0 & 0\\ 0 & \boldsymbol{1}\end{bmatrix}.
\]
\hfill\qedsymbol
\end{example}
As in Example \ref{Ex1}, we observe that taking $(\varPsi_t\colon
t\geq 0)$ --- an ergodic semigroup of positive maps on $M$ with
$\sigma((\varPsi_t)) = \{1,-1\}$, and defining
\[
 \varPhi_t\left(\mat\right) =
 \begin{bmatrix}\varPsi_t\left(\frac{a+d}{2}\right) &
 \lambda_0^t\varPsi_t(b) \\ \bar{\lambda}_0^t\varPsi_t(c) &
 \varPsi_t\left(\frac{a+d}{2}\right)\end{bmatrix},\qquad t\geq 0,
\]
we obtain a continuous counterpart of Example \ref{Ex2}.

\end{document}